\documentclass[11pt,english]{article}

\usepackage[latin9]{inputenc}
\usepackage{amsmath}
\usepackage{amssymb}
\usepackage[a4paper]{geometry}
\usepackage{setspace}
\PassOptionsToPackage{normalem}{ulem}
\usepackage{ulem}
\makeatletter

\providecommand{\tabularnewline}{\\}

\@ifundefined{date}{}{\date{}}
\usepackage{array}
\usepackage{mathrsfs}
\usepackage{color}

\makeatother

\usepackage{babel}
\begin{document}
\title{\textbf{On the Algebraic Independence of a}\\
\textbf{Set of Generalized Constants}}
\author{Michael R. Powers\thanks{Department of Finance, School of Economics and Management, and Schwarzman
College, Tsinghua University, Beijing, China 100084; email: powers@sem.tsinghua.edu.cn.}}
\date{August 23, 2025}
\maketitle
\begin{abstract}
\begin{singlespace}
\noindent Neither the Euler-Mascheroni constant, $\gamma=0.577215\ldots$,
nor the Euler-Gompertz constant, $\delta=0.596347\ldots$, is currently
known to be irrational. However, it has been proved that these two
numbers are disjunctively transcendental; that is, at least one of
them must be transcendental. The two constants are related through
a well-known equation of Hardy, which recently has been generalized
to a pair of infinite sequences $(\gamma^{\left(n\right)},\delta^{\left(n\right)})$
based on moments of the Gumbel(0,1) probability distribution. In the
present work, we demonstrate the algebraic independence of the set
$\{\gamma^{\left(n\right)}+\delta^{\left(n\right)}/e\}_{n\geq0}$,
and thus the transcendence of $\gamma^{\left(n\right)}+\delta^{\left(n\right)}/e$
for all $n\geq0$. This further implies the disjunctive transcendence
of both pairs $(\gamma^{\left(n\right)},\delta^{\left(n\right)}/e)$
and $(\gamma^{\left(n\right)},\delta^{\left(n\right)})$ for all $n\geq1$.\medskip{}

\end{singlespace}
\end{abstract}
\begin{singlespace}

\section{Introduction}
\end{singlespace}

\begin{singlespace}
\noindent Let $\gamma=0.577215\ldots$ denote the Euler-Mascheroni
constant and $\delta=0.596347\ldots$ denote the Euler-Gompertz constant.
Although neither $\gamma$ nor $\delta$ has been shown to be irrational,
Shidlovski (1962) implicitly proved their disjunctive irrationality;\footnote{\noindent This disjunctive-irrationality result is sometimes attributed
to Mahler (1968) (who was the first to state it explicitly) and/or
Aptekarev (2009) (who quantified it within a linear system). } that is, at least one of the two numbers must be irrational. More
recently, Rivoal (2012) strengthened this result to show the disjunctive
transcendence of $\left(\gamma,\delta\right)$.\footnote{\noindent See Lagarias (2013) for a comprehensive treatment of results
involving these two constants.}
\end{singlespace}

\begin{singlespace}
The two constants are related by the intriguing equation
\begin{equation}
{\displaystyle \sum_{k=1}^{\infty}\dfrac{\left(-1\right)^{k+1}}{k\cdot k!}}=\gamma+\dfrac{\delta}{e},
\end{equation}
first noted by Hardy (1949). Given that $\delta=-e\textrm{Ei}\left(-1\right)$,
this expression follows immediately by setting $u=-1$ in the power-series
expansion of the exponential-integral function,
\[
\textrm{Ei}\left(u\right)=\gamma+\ln\left|u\right|+{\displaystyle \sum_{k=1}^{\infty}\dfrac{u^{k}}{k\cdot k!}}.
\]
Recently, Powers (2025) provided the following probabilistic interpretation
of (1):
\begin{equation}
\textrm{E}_{X}\left[X^{+}\right]=\textrm{E}_{X}\left[X\right]-\textrm{E}_{X}\left[X\mid X\leq0\right]\Pr\left\{ X\leq0\right\} ,
\end{equation}
where $X\sim\textrm{Gumbel}\left(0,1\right)$\footnote{\noindent The cumulative distribution function of $\textrm{Gumbel}\left(0,1\right)$
is given by $F_{X}\left(x\right)=\exp\left(-e^{-x}\right),\:x\in\mathbb{R}$,
with mean $\gamma$ and variance $\pi^{2}/6$.} so that
\[
\textrm{E}_{X}\left[X^{+}\right]={\displaystyle \int}_{0}^{\infty}x\exp\left(-x-e^{-x}\right)dx={\displaystyle \sum_{k=1}^{\infty}\dfrac{\left(-1\right)^{k+1}}{k\cdot k!}},
\]
\[
\textrm{E}_{X}\left[X\right]={\displaystyle \int}_{-\infty}^{\infty}x\exp\left(-x-e^{-x}\right)dx=\gamma,
\]
\[
-\textrm{E}_{X}\left[X\mid X\leq0\right]=-e{\displaystyle \int_{-\infty}^{0}}x\exp\left(-x-e^{-x}\right)dx=\delta,
\]
and
\[
\Pr\left\{ X\leq0\right\} ={\displaystyle \int_{-\infty}^{0}}\exp\left(-x-e^{-x}\right)dx=\dfrac{1}{e}.
\]
The identity in (2) suggests the natural generalization
\[
\textrm{E}_{X}\left[\left(X^{+}\right)^{n}\right]=\textrm{E}_{X}\left[X^{n}\right]-\textrm{E}_{X}\left[X^{n}\mid X\leq0\right]\Pr\left\{ X\leq0\right\} 
\]
for $n\in\mathbb{Z}_{\geq0}$, which can be expressed in the manner
of (1) as
\begin{equation}
n!{\displaystyle \sum_{k=1}^{\infty}\dfrac{\left(-1\right)^{k+1}}{k^{n}\cdot k!}}=\gamma^{\left(n\right)}+\dfrac{\delta^{\left(n\right)}}{e},
\end{equation}
where
\[
n!{\displaystyle \sum_{k=1}^{\infty}\dfrac{\left(-1\right)^{k+1}}{k^{n}\cdot k!}}=\textrm{E}_{X}\left[\left(X^{+}\right)^{n}\right]={\displaystyle \int}_{0}^{\infty}x^{n}\exp\left(-x-e^{-x}\right)dx,
\]
\[
\gamma^{\left(n\right)}=\textrm{E}_{X}\left[X^{n}\right]={\displaystyle \int}_{-\infty}^{\infty}x^{n}\exp\left(-x-e^{-x}\right)dx,
\]
and
\begin{equation}
\delta^{\left(n\right)}=-\textrm{E}_{X}\left[X^{n}\mid X\leq0\right]=-e{\displaystyle \int_{-\infty}^{0}}x^{n}\exp\left(-x-e^{-x}\right)dx.
\end{equation}

For ease of exposition, we will set
\begin{equation}
\eta^{\left(n\right)}=n!{\displaystyle \sum_{k=1}^{\infty}\dfrac{\left(-1\right)^{k+1}}{k^{n}\cdot k!}}
\end{equation}
so that (3) may be rewritten as
\begin{equation}
\eta^{\left(n\right)}=\gamma^{\left(n\right)}+\dfrac{\delta^{\left(n\right)}}{e}.
\end{equation}
Further, we will refer to the $\eta^{\left(n\right)}$, $\gamma^{\left(n\right)}$
, and $\delta^{\left(n\right)}$, as \emph{generalized Eta}, \emph{Euler-Mascheroni},
and \emph{Euler-Gompertz} constants, respectively (with $\eta^{\left(1\right)}=\eta$,
$\gamma^{\left(1\right)}=\gamma$ , and $\delta^{\left(1\right)}=\delta$
denoting the \emph{ordinary} Eta, Euler-Mascheroni, and Euler-Gompertz
constants).\footnote{Note that the generalized Euler-Mascheroni constants, $\gamma^{\left(n\right)}$,
are indexed via superscripts to distinguish them from the well-known
sequence of Stieltjes constants, commonly denoted by $\gamma_{n}$.} The following section summarizes certain basic properties of the
three sequences. Then, in Section 3, we prove: (i) the set of numbers
$\left\{ \eta^{\left(n\right)}\right\} _{n\geq0}$ are algebraically
independent, implying that $\eta^{\left(n\right)}$ is transcendental
for all $n\in\mathbb{Z}_{\geq0}$; and (ii) both constant pairs $\left(\gamma^{\left(n\right)},\delta^{\left(n\right)}/e\right)$
and $(\gamma^{\left(n\right)},\delta^{\left(n\right)})$ are disjunctively
transcendental for all $n\in\mathbb{Z}_{\geq1}$.\medskip{}
\medskip{}

\end{singlespace}
\begin{singlespace}

\section{The Sequences $\boldsymbol{\eta^{\left(n\right)}}$, $\boldsymbol{\gamma^{\left(n\right)}}$,
and $\boldsymbol{\delta^{\left(n\right)}}$}
\end{singlespace}

\begin{singlespace}
\noindent Table 1 presents values of the first 16 generalized Euler-Mascheroni,
Euler-Gompertz, and Eta constants.
\end{singlespace}
\begin{singlespace}
\begin{center}
Table 1. Values of $\eta^{\left(n\right)}$, $\gamma^{\left(n\right)}$,
and $\delta^{\left(n\right)}$ for $n\in\left\{ 0,1,\ldots,15\right\} $
\par\end{center}

\begin{center}
\begin{tabular}{|c|c|c|c|}
\hline 
{\small$\boldsymbol{n}$} & {\small$\boldsymbol{\eta^{\left(n\right)}}$} & {\small$\boldsymbol{\gamma^{\left(n\right)}}$} & {\small$\boldsymbol{\delta^{\left(n\right)}}$}\tabularnewline
\hline 
\hline 
{\small$0$} & {\small$0.6321205588\ldots$} & {\small$1.0$} & {\small$-1.0$}\tabularnewline
\hline 
{\small$1$} & {\small$0.7965995992\ldots$} & {\small$0.5772156649\ldots$} & {\small$0.5963473623\ldots$}\tabularnewline
\hline 
{\small$2$} & {\small$1.7824255962\ldots$} & {\small$1.9781119906\ldots$} & {\small$\text{\textendash}0.5319307700\ldots$}\tabularnewline
\hline 
{\small$3$} & {\small$5.6584954080\ldots$} & {\small$5.4448744564\ldots$} & {\small$0.5806819508\ldots$}\tabularnewline
\hline 
{\small$4$} & {\small$23.2957725933\ldots$} & {\small$23.5614740840\ldots$} & {\small$\text{\textendash}0.7222515339\ldots$}\tabularnewline
\hline 
{\small$5$} & {\small$118.2027216118\ldots$} & {\small$117.8394082683\ldots$} & {\small$0.9875880596\ldots$}\tabularnewline
\hline 
{\small$6$} & {\small$714.5326485509\ldots$} & {\small$715.0673625273\ldots$} & {\small$\text{\textendash}1.4535032853\ldots$}\tabularnewline
\hline 
{\small$7$} & {\small$5,020.6842841603\ldots$} & {\small$5,019.8488726298\ldots$} & {\small$2.2708839827\ldots$}\tabularnewline
\hline 
{\small$8$} & {\small$40,242.2494274946\ldots$} & {\small$40,243.6215733357\ldots$} & {\small$\text{\textendash}3.7298791058\ldots$}\tabularnewline
\hline 
{\small$9$} & {\small$362,528.6415241055\ldots$} & {\small$362,526.2891146549\ldots$} & {\small$6.3945118625\ldots$}\tabularnewline
\hline 
{\small$10$} & {\small$\qquad3,627,038.2261612415\ldots\qquad$} & {\small$\qquad3,627,042.4127568947\ldots\qquad$} & {\small$\text{\textendash}11.3803468877\ldots$}\tabularnewline
\hline 
{\small$11$} & {\small$39,907,091.8528764918\ldots$} & {\small$39,907,084.1514313358\ldots$} & {\small$20.9346984188\ldots$}\tabularnewline
\hline 
{\small$12$} & {\small$478,943,277.1724405288\ldots$} & {\small$478,943,291.7651829432\ldots$} & {\small$\text{\textendash}39.6671864816\ldots$}\tabularnewline
\hline 
{\small$13$} & {\small$6,226,641,379.9457960128\ldots$} & {\small$6,226,641,351.5460642549\ldots$} & {\small$77.1984745660\ldots$}\tabularnewline
\hline 
{\small$14$} & {\small$87,175,633,754.0756530761\ldots$} & {\small$87,175,633,810.7084156319\ldots$} & {\small$\text{\textendash}153.9437943882\ldots$}\tabularnewline
\hline 
{\small$15$} & {\small$1,307,654,429,611.2775878906\ldots$} & {\small$1,307,654,429,495.7941762096\ldots$} & {\small$313.9164765016\ldots$}\tabularnewline
\hline 
\end{tabular}
\par\end{center}
\end{singlespace}
\begin{singlespace}

\subsection{Asymptotic Behavior}
\end{singlespace}

\begin{singlespace}
From the columns of Table 1, it is clear that both the $\eta^{\left(n\right)}$
and $\gamma^{\left(n\right)}$ increase approximately factorially,
whereas the absolute magnitudes of the $\delta^{\left(n\right)}$,
which alternate in sign, appear to grow slightly super-exponentially.
The following proposition refines these observations.\medskip{}
\medskip{}

\end{singlespace}

\begin{singlespace}
\noindent\textbf{Proposition 1:} As $n\rightarrow\infty$:
\end{singlespace}

\begin{singlespace}
\[
\left(\textrm{i}\right)\:\eta^{\left(n\right)}=n!\left(1-\dfrac{1}{2^{n+1}}+O\left(\dfrac{1}{3^{n}}\right)\right);
\]

\[
\left(\textrm{ii}\right)\:\gamma^{\left(n\right)}=n!\left(1-\dfrac{1}{2^{n+1}}+O\left(\dfrac{1}{3^{n}}\right)\right);\textrm{and}
\]

\[
\left(\textrm{iii}\right)\:\left|\delta^{\left(n\right)}\right|=\left(\ln\left(n\right)\right)^{n}\exp\left(O\left(\dfrac{n}{\ln\left(n\right)}\right)\right).
\]

\medskip{}
\medskip{}

\end{singlespace}

\begin{singlespace}
\noindent\textbf{Proof:} We consider the three parts of the proposition
in the order: (i), (iii), (ii).\medskip{}

\noindent (i) Consider the right-hand side of (5):
\[
\eta^{\left(n\right)}=n!\left(1-\dfrac{1}{2^{n}\cdot2!}+\dfrac{1}{3^{n}\cdot3!}-\dfrac{1}{4^{n}\cdot4!}+\cdots\right).
\]
From this expansion it is easy to see that
\[
\eta^{\left(n\right)}<n!\left[1-\dfrac{1}{2^{n+1}}+\left(\dfrac{1}{6}\right)\dfrac{1}{3^{n}}+\left(\dfrac{1}{24}\right)\dfrac{1}{4^{n}}+\cdots\right]
\]
\[
<n!\left[1-\dfrac{1}{2^{n+1}}+\left(\dfrac{1}{6}+\dfrac{1}{24}+\cdots\right)\dfrac{1}{3^{n}}\right]
\]
\[
=n!\left[1-\dfrac{1}{2^{n+1}}+\dfrac{\left(e-5/2\right)}{3^{n}}\right]
\]
and
\[
\eta^{\left(n\right)}>n!\left[1-\dfrac{1}{2^{n+1}}-\left(\dfrac{1}{6}\right)\dfrac{1}{3^{n}}-\left(\dfrac{1}{24}\right)\dfrac{1}{4^{n}}-\cdots\right]
\]
\[
>n!\left[1-\dfrac{1}{2^{n+1}}-\left(\dfrac{1}{6}+\dfrac{1}{24}+\cdots\right)\dfrac{1}{3^{n}}\right]
\]
\[
=n!\left[1-\dfrac{1}{2^{n+1}}-\dfrac{\left(e-5/2\right)}{3^{n}}\right].
\]
It then follows that
\[
\eta^{\left(n\right)}=n!\left(1-\dfrac{1}{2^{n+1}}+O\left(\dfrac{1}{3^{n}}\right)\right).
\]
\medskip{}

\noindent (iii) Substituting $u=e^{-x}$ into the integral in (4)
yields
\[
\delta^{\left(n\right)}=\left(-1\right)^{n+1}e{\displaystyle {\displaystyle \int_{1}^{\infty}}\left[\ln\left(u\right)\right]^{n}e^{-u}du}
\]
\[
\Longrightarrow\left|\delta^{\left(n\right)}\right|=e{\displaystyle {\displaystyle \int_{1}^{\infty}}\left[\ln\left(u\right)\right]^{n}e^{-u}du}
\]
\begin{equation}
=e{\displaystyle {\displaystyle \int_{1}^{\infty}}e^{\phi_{n}\left(u\right)}du},
\end{equation}
where $\phi_{n}\left(u\right)=n\ln\left(\ln\left(u\right)\right)-u$.
We then employ Laplace's (saddle-point) method to approximate (7).
\end{singlespace}

\begin{singlespace}
Taking derivatives of $\phi_{n}\left(u\right)$ with respect to $u$
gives
\[
\phi_{n}^{\prime}\left(u\right)=\dfrac{n}{u\ln\left(u\right)}-1
\]
and
\[
\phi_{n}^{\prime\prime}\left(u\right)=-\dfrac{n}{u^{2}\ln\left(u\right)}\left(1+\dfrac{1}{\ln\left(u\right)}\right)<0,
\]
 revealing that $\phi_{n}\left(u\right)$ enjoys a unique global maximum
at the saddle point
\[
u^{*}\ln\left(u^{*}\right)=n
\]
\[
\Longleftrightarrow u^{*}=\dfrac{n}{W\left(n\right)},\:\ln\left(u^{*}\right)=W\left(n\right).
\]
Then
\[
\phi_{n}\left(u^{*}\right)=n\ln\left(\ln\left(u^{*}\right)\right)-\dfrac{n}{\ln\left(u^{*}\right)}
\]
\[
=n\ln\left(W\left(n\right)\right)-\dfrac{n}{W\left(n\right)},
\]
\[
\left|\phi_{n}^{\prime\prime}\left(u^{*}\right)\right|=\dfrac{W\left(n\right)}{n}\left(1+\dfrac{1}{W\left(n\right)}\right)
\]
\[
=\dfrac{W\left(n\right)+1}{n},
\]
and (7) can be approximated by
\[
\left|\delta^{\left(n\right)}\right|=e\cdot e^{\phi_{n}\left(u^{*}\right)}\sqrt{\dfrac{2\pi}{\left|\phi_{n}^{\prime\prime}\left(u^{*}\right)\right|}}\left(1+o\left(1\right)\right)
\]
\[
=e\left[W\left(n\right)\right]^{n}\exp\left(-\dfrac{n}{W\left(n\right)}\right)\sqrt{\dfrac{2\pi n}{W\left(n\right)+1}}\left(1+o\left(1\right)\right)
\]
\[
=e\left(\ln\left(n\right)\right)^{\left(1+o\left(1\right)\right)n}\exp\left(-\dfrac{n}{\ln\left(n\right)}\left(1+o\left(1\right)\right)\right)\sqrt{\dfrac{2\pi n}{\ln\left(n\right)}}\left(1+O\left(\dfrac{\ln\left(\ln\left(n\right)\right)}{\ln\left(n\right)}\right)\right)\left(1+o\left(1\right)\right)
\]
\[
=\left(\ln\left(n\right)\right)^{n}\exp\left(O\left(\dfrac{n}{\ln\left(n\right)}\right)\right).
\]
\medskip{}

\end{singlespace}

\begin{singlespace}
\noindent (ii) Finally, we assemble the results in (i) and (iii) via
(6), giving
\[
\gamma^{\left(n\right)}=n!\left(1-\dfrac{1}{2^{n+1}}+O\left(\dfrac{1}{3^{n}}\right)\right)-\left(-1\right)^{n+1}\left(\ln\left(n\right)\right)^{n}\exp\left(O\left(\dfrac{n}{\ln\left(n\right)}\right)\right)
\]
\[
=n!\left(1-\dfrac{1}{2^{n+1}}+O\left(\dfrac{1}{3^{n}}\right)\right)+o\left(n!\right)
\]
\[
=n!\left(1-\dfrac{1}{2^{n+1}}+O\left(\dfrac{1}{3^{n}}\right)\right).\:\blacksquare
\]

\end{singlespace}
\begin{singlespace}

\subsection{A Recurrence Relation}
\end{singlespace}

\begin{singlespace}
\noindent Consider the recurrence
\[
\textrm{E}_{X}\left[X^{n+1}\right]={\displaystyle \sum_{j=0}^{n}}\binom{n}{j}\kappa_{n+1-j}\textrm{E}_{X}\left[X^{j}\right],
\]
derived from Bell polynomials, which holds for any random variable
$X$ with finite $\left(n+1\right)^{\textrm{st}}$ raw moment (where
$\kappa_{j}$ denotes the $j^{\textrm{th}}$ cumulant). This identity
can be tailored for $X\sim\textrm{Gumbel}\left(0,1\right)$ by inserting
\[
\kappa_{j}=\begin{cases}
\gamma & j=1\\
\left(j-1\right)!\zeta\left(j\right), & j\in\left\{ 2,3,\ldots\right\} 
\end{cases},
\]
immediately yielding
\begin{equation}
\gamma^{\left(n+1\right)}=\gamma\cdot\gamma^{\left(n\right)}+{\displaystyle \sum_{j=0}^{n-1}}\dfrac{n!}{j!}\zeta\left(n+1-j\right)\gamma^{\left(j\right)}.
\end{equation}
Then, (8) allows us to solve for the individual generalized Euler-Mascheroni
constants as follows:
\[
\gamma^{\left(1\right)}=\gamma,
\]
\[
\gamma^{\left(2\right)}=\gamma^{2}+\zeta\left(2\right),
\]
\[
\gamma^{\left(3\right)}=\gamma^{3}+3\gamma\zeta\left(2\right)+2\zeta\left(3\right),
\]
\[
\gamma^{\left(4\right)}=\gamma^{4}+6\zeta\left(2\right)\gamma^{2}+8\zeta\left(3\right)\gamma+\dfrac{27}{2}\zeta\left(4\right),
\]
\[
\textrm{etc.}
\]

\end{singlespace}

\begin{singlespace}
Unfortunately, neither the generalized Eta nor generalized Euler-Gompertz
constants appear to satisfy simple recurrences such as (8).\medskip{}
\medskip{}

\end{singlespace}
\begin{singlespace}

\section{Main Results}
\end{singlespace}

\begin{singlespace}
\noindent The core results of our study follow from the sequence of
functions
\begin{equation}
F_{n}\left(t\right)={\displaystyle \int_{0}^{\infty}}x^{n}\exp\left(-\left(x+te^{-x}\right)\right)dx,
\end{equation}
for $n\in\mathbb{Z}_{\geq0}$ and $t\in\mathbb{R}^{+}$. These may
be shown to be $E$-functions (in the sense of Siegel) because they
satisfy the four properties below.\medskip{}

\noindent (1) \uline{Entirety}:
\[
F_{n}\left(t\right)={\displaystyle \int_{0}^{\infty}}x^{n}e^{-x}{\displaystyle \sum_{k=0}^{\infty}}\dfrac{\left(-t\right)^{k}e^{-kx}}{k!}dx
\]
\[
={\displaystyle \sum_{k=0}^{\infty}}\dfrac{\left(-t\right)^{k}}{k!}{\displaystyle \int_{0}^{\infty}}x^{n}e^{-\left(k+1\right)x}dx
\]
(where the integral and summation can be interchanged by the Dominated
Convergence theorem because $\left|x^{n}e^{-x}{\textstyle \sum_{k=0}^{\infty}\left(-t\right)^{k}e^{-kx}/k!}\right|\leq x^{n}e^{-x}\exp\left(\left|t\right|e^{-x}\right)\leq x^{n}e^{-x+\left|t\right|}$
for all $x$)
\[
={\displaystyle \sum_{k=0}^{\infty}}\dfrac{n!\left(-t\right)^{k}}{\left(k+1\right)^{n+1}k!}
\]
\[
={\displaystyle \sum_{k=0}^{\infty}}\dfrac{n!\left(-t\right)^{k}}{\left(k+1\right)^{n}\left(k+1\right)!}
\]
\[
={\displaystyle \sum_{k=0}^{\infty}}\dfrac{a_{k}t^{k}}{k!}
\]
holds for all $t\in\mathbb{C}$.\medskip{}

\noindent (2) \uline{Bound on Algebraic Conjugates}:
\[
a_{k}=\dfrac{\left(-1\right)^{k}n!}{\left(k+1\right)^{n+1}}\in\mathbb{Q}
\]
for all $k\in\mathbb{Z}_{\geq0}$, meaning that the algebraic conjugate
of $a_{k}$ is $a_{k}$ itself. We then see $\left|a_{k}\right|\leq n!$
for all $k$, which implies there exists a constant $C>0$ such that
$\left|a_{k}\right|\leq C^{k+1}$ for all $k$.\medskip{}

\noindent (3) \uline{Polynomial-Coefficient Linear Differential Equation}:
\[
\widetilde{D}\left[F_{n}\left(t\right)\right]=\begin{cases}
e^{-t}, & n=0\\
nF_{n-1}\left(t\right), & n\geq1
\end{cases}
\]
for $\widetilde{D}\equiv t\left(d/dt\right)+1$ implies
\begin{equation}
\widetilde{D}^{n+1}\left[F_{n}\left(t\right)\right]=n!e^{-t}.
\end{equation}
Furthermore,
\begin{equation}
D\left[n!e^{-t}\right]=0
\end{equation}
for $D\equiv\left(d/dt\right)+1$. Combining (10) and (11) then yields
the linear differential equation
\[
D\widetilde{D}^{n+1}\left[F_{n}\left(t\right)\right]=0
\]
with polynomial coefficients in $\overline{\mathbb{Q}}\left[t\right]$.\medskip{}

\noindent (4) \uline{Bound on Denominator Growth}:
\[
g_{k}=\left(\textrm{lcm}\left(1,2,\ldots,k+1\right)\right)^{n+1}
\]
is a common denominator of $a_{0},a_{1},\ldots,a_{k}$ for $k\in\mathbb{Z}_{\geq0}$.
Since
\[
\ln\left(\textrm{lcm\ensuremath{\left(1,2,\ldots,k+1\right)}}\right)=\psi\left(k+1\right)
\]
\[
=k+1+o\left(k\right)
\]
(where $\psi\left(\cdot\right)$ denotes Chebyshev's second function),
it follows that
\[
\ln\left(g_{k}\right)=\left(n+1\right)\left(k+1\right)+o\left(k\right).
\]
This implies there exists a constant $C>0$ such that
\[
g_{k}\leq\left[e^{C\left(n+1\right)}\right]^{k+1}
\]
for all $k$.\medskip{}

\end{singlespace}

\begin{singlespace}
Now consider the following lemma, necessary for proving Theorems 1
and 2 below.

\medskip{}

\end{singlespace}

\begin{singlespace}
\noindent\textbf{Lemma 1:} For all $m\in\mathbb{Z}_{\geq1}$, the
$E$-functions $F_{0}\left(t\right),F_{1}\left(t\right),\ldots,F_{m}\left(t\right)$
are algebraically independent over $\mathbb{C}\left(t\right)$.\medskip{}
\medskip{}

\noindent\textbf{Proof:} We construct a sequence of differential
field extensions, each equipped with the same derivation, $T\equiv t\left(d/dt\right)$.
\end{singlespace}

\begin{singlespace}
Let $\mathcal{K}=\mathbb{C}\left(t\right)$ be the base differential
field (of characteristic 0) with algebraically closed field of constants
$\mathbb{C}$, and consider the system of $E$-functions given by
(9). First, note that
\begin{equation}
F_{0}^{\prime}\left(t\right)=\dfrac{1}{t}-\left(1+\dfrac{1}{t}\right)F_{0}\left(t\right),
\end{equation}
where $F_{0}\left(t\right)=\left(1-e^{-t}\right)/t$, and
\begin{equation}
F_{n}^{\prime}\left(t\right)=\dfrac{1}{t}\left(nF_{n-1}\left(t\right)-F_{n}\left(t\right)\right)
\end{equation}
for $n\geq1$. Next, set
\[
E\left(t\right)=1-tF_{0}\left(t\right)
\]
\begin{equation}
=e^{-t}
\end{equation}
and
\begin{equation}
H_{n}\left(t\right)=tF_{n}\left(t\right)
\end{equation}
for $n\in\mathbb{Z}_{\geq0}$, so that
\[
E^{\prime}\left(t\right)=-E\left(t\right),
\]
\[
H_{0}^{\prime}\left(t\right)=E\left(t\right),
\]
and
\[
H_{n}^{\prime}\left(t\right)=\dfrac{n}{t}H_{n-1}\left(t\right)
\]
for $n\geq1$. Then
\[
T\left[E\right]=-tE
\]
and
\[
T\left[H_{n}\right]=nH_{n-1}
\]
for $n\geq1$.

Now note that if $L$ is any differential field extension of $\mathcal{K}$,
then
\[
\textrm{Const}\left(L,T\right)=\left\{ y\in L:T\left[y\right]=0\right\} 
\]
\[
=\left\{ y\in L:t\dfrac{d}{dt}\left[y\right]=0\right\} 
\]
\[
=\left\{ y\in L:\dfrac{d}{dt}\left[y\right]=0\right\} 
\]
\[
=\textrm{Const}\left(L,\dfrac{d}{dt}\right)
\]
because $t\in\mathcal{K}^{\times}\subset L^{\times}$. Thus, $T$
and $d/dt$ have the same constants on every field we adjoin. In particular,
\[
\textrm{Const}\textrm{\ensuremath{\left(\mathcal{K},T\right)}}=\textrm{Const}\textrm{\ensuremath{\left(\mathcal{K},d/dt\right)}}
\]
\[
=\mathbb{C}.
\]

The following four Facts are necessary to complete the proof.\medskip{}
\medskip{}

\end{singlespace}

\begin{singlespace}
\noindent\textbf{Fact 1:} $E$ is transcendental over $\mathcal{K}$.

\noindent Since $\left(d/dt\right)\left[E\right]/E=-1\in\mathcal{K}$,
it follows that $E$ is an exponential of an integral over $\mathcal{K}$.
If such an element were algebraic over the base (with characteristic
0 and algebraically closed constants), then the Kolchin-Ostrowski
theorem for exponentials of integrals (e.g., Theorem 2.6 of Srinivasan,
2008) would imply that $E^{n}\in\mathcal{K}$ for some $n\geq1$,
which is impossible. (Note that the Kolchin-Ostrowski theorem is formulated
for an arbitrary derivation, and we use $d/dt$.) Therefore, $E$
must be transcendental. $\blacklozenge$\medskip{}
\medskip{}

\noindent\textbf{Fact 2:} There is no enlargement of constants along
the tower.

\noindent Define the tower (with fixed derivation $T$):
\[
\mathcal{K}_{0}=\mathcal{K},
\]
\[
\mathcal{K}_{1}=\mathcal{K}\left(E\right),
\]
and
\[
\mathcal{K}_{n+1}=\mathcal{K}_{n}\left(H_{n}\right)
\]
for $n\geq1$. Given that $T\left[E\right]/E=-t\in\mathcal{K}_{0}$,
the Kolchin-Ostrowski theorem for exponentials of integrals (e.g.,
Theorem 2.6 of Srinivasan, 2008) implies that adjoining $E$ to $\mathcal{K}_{0}$
(with characteristic 0 and algebraically closed constants) does not
adjoin new constants. (Again, recall that the Kolchin-Ostrowski theorem
holds for an arbitrary derivation, in particular $T$.) Thus,
\begin{equation}
\textrm{Const}\textrm{\ensuremath{\left(\mathcal{K}_{1},T\right)}}=\mathbb{C}.
\end{equation}
For $n\geq1$, we prove the following two Assertions by induction:
\[
\left(\textrm{I}_{n}\right)\:\textrm{Const}\left(\mathcal{K}_{n},T\right)=\mathbb{C};
\]
and
\[
\left(\textrm{II}_{n}\right)\:\textrm{There is no }G\in\mathcal{K}_{n}\textrm{ such that }T\left[G\right]=nH_{n-1}.
\]

\end{singlespace}

\begin{singlespace}
For the base case ($n=1$), Assertion $\left(\textrm{I}_{1}\right)$
holds by (16) and Assertion $\left(\textrm{II}_{1}\right)$ is equivalent
to Fact 3 (proved below; that is, there is no $G\in\mathcal{K}_{1}=\mathcal{K}(E)$
such that $T\left[G\right]=H_{0}=1-E$). Now assume that both Assertions
$\left(\textrm{I}_{n}\right)$ and $\left(\textrm{II}_{n}\right)$
hold for some fixed $n\geq1$, and set $\widetilde{\mathcal{K}}=\mathcal{K}_{n}$
and $A=H_{n}$ so that $T\left[A\right]=nH_{n-1}\in\widetilde{\mathcal{K}},$
$\mathrm{Const}\left(\widetilde{\mathcal{K}},T\right)=\mathbb{C},$
and $nH_{n-1}\notin T\left[\widetilde{\mathcal{K}}\right]$. We will
proceed in three steps:\medskip{}

\end{singlespace}

\begin{singlespace}
\noindent\uline{Step 1}: No enlargement of constants when adjoining
$A$.

\noindent Let $L=\widetilde{\mathcal{K}}\left(A\right)$ and suppose
$B\in L$ with $T\left[B\right]=0$. We now show $B\in\mathbb{C}$.
\end{singlespace}

\begin{singlespace}
\emph{Case 1: $A$ is transcendental over $\widetilde{\mathcal{K}}$.}
Write $B=P\left(A\right)/Q\left(A\right)$ with co-prime $P,Q\in\widetilde{\mathcal{K}}\left[A\right]$.
If $\deg_{A}\left(Q\right)>0$, then a standard valuation/pole-order
test shows that if $A$ is transcendental over a differential field
$L$ and $T\left[A\right]\in L$, then $T\left[B\right]$ has a pole
at a finite point (because differentiation of a proper rational function
in $A$ increases the pole order), contradicting $T\left[B\right]=0$.
Thus, $Q\in\widetilde{\mathcal{K}}^{\times}$ and
\[
B=\sum_{j=0}^{d}b_{j}A^{j}
\]
with $b_{j}\in\widetilde{\mathcal{K}}$, giving
\[
T\left[B\right]=\sum_{j=0}^{d}T\left[b_{j}\right]A^{j}+\sum_{j=1}^{d}jb_{j}\left(nH_{n-1}\right)A^{j-1}
\]
\[
=0.
\]
Comparing coefficients at the highest $A$ degree reveals $T\left[b_{d}\right]=0$,
so $b_{d}\in\mathbb{C}$ by Assertion $\left(\textrm{I}_{n}\right)$.
Then, the coefficients at $A^{d-1}$ show that $T\left[b_{d-1}\right]+db_{d}\left(nH_{n-1}\right)=0$,
forcing $b_{d}=0$ because $nH_{n-1}\notin T\left[\widetilde{\mathcal{K}}\right]$
by Assertion $\left(\textrm{II}_{n}\right)$. Descending in degree,
we obtain $b_{j}=0$ for all $j\ge1$, and finally $T\left[b_{0}\right]=0$,
implying $b_{0}\in\mathbb{C}$. Therefore, $B\in\mathbb{C}$.\medskip{}

\emph{Case 2: $A$ is algebraic over $\widetilde{\mathcal{K}}$.}
Let $d=\dim_{\widetilde{\mathcal{K}}}\left(\widetilde{\mathcal{K}}\left(A\right)\right)$,
and write
\[
B=\sum_{j=0}^{d-1}b_{j}A^{j},
\]
the $\widetilde{\mathcal{K}}$-basis $\left\{ 1,A,A^{2},\dots,A^{d-1}\right\} $.
Then
\[
T\left[B\right]=\sum_{j=0}^{d-1}T\left[b_{j}\right]A^{j}+\sum_{j=1}^{d-1}jb_{j}\left(nH_{n-1}\right)A^{j-1}
\]
\[
=0,
\]
and by the $\widetilde{\mathcal{K}}$-linear independence of $\left\{ 1,A,A^{2},\dots,A^{d-1}\right\} $,
the coefficients at $A^{d-1}$ yield $T\left[b_{d-1}\right]=0$, so
$b_{d-1}\in\mathbb{C}$ by Assertion $\left(\textrm{I}_{n}\right)$.
Further, the coefficients at $A^{d-2}$ give $T\left[b_{d-2}\right]+\left(d-1\right)b_{d-1}\left(nH_{n-1}\right)=0$,
forcing $b_{d-1}=0$ because $nH_{n-1}\notin T\left[\widetilde{\mathcal{K}}\right]$
by Assertion $\left(\textrm{II}_{n}\right)$. Proceeding downward,
we conclude $b_{j}=0$ for all $j\ge1$, and finally $T\left[b_{0}\right]=0$,
so that $b_{0}\in\mathbb{C}$. Therefore, $B\in\mathbb{C}$.

In either of the above Cases, we have shown
\[
\mathrm{Const}\left(L,T\right)=\mathbb{C},
\]
thus proving Assertion $\left(\textrm{I}_{n+1}\right)$.\medskip{}

\end{singlespace}

\begin{singlespace}
\noindent\uline{Step 2}: $A=H_{n}$ is transcendental over $\widetilde{\mathcal{K}}$\emph{.}

\noindent If $A$ were algebraic over $\widetilde{\mathcal{K}}$,
then $L/\widetilde{\mathcal{K}}$ would be a finite algebraic differential
extension with the same constants $\mathbb{C}$ (by Step 1). Applying
the standard trace argument for finite algebraic differential extensions
with equal constants to $L=\mathcal{K}_{n+1}$, $\widetilde{\mathcal{K}}=\mathcal{K}_{n}$,
and the element $A=H_{n}\in L$ (for which $T\left[A\right]=nH_{n-1}\in\widetilde{\mathcal{K}}$)
then would give $A\in\widetilde{\mathcal{K}}$, contradicting Assertion
$\left(\textrm{II}_{n}\right)$. Therefore, $A$ is transcendental
over $\widetilde{\mathcal{K}}$.\medskip{}
\medskip{}

\noindent\uline{Step 3}: Assertion $\left(\textrm{II}_{n+1}\right)$
is valid.

\noindent Assume (for purposes of contradiction) that there exists
$G\in L$ with $T\left[G\right]=\left(n+1\right)A$, where $L=\widetilde{\mathcal{K}}\left(A\right)$
and $A=H_{n}$. Since $A$ is transcendental over $\widetilde{\mathcal{K}}$
(by Step 2) and $T\left[A\right]=nH_{n-1}\in\widetilde{\mathcal{K}}$,
we can write $G=P\left(A\right)/Q\left(A\right)$ with co-prime $P,Q\in\widetilde{\mathcal{K}}\left[A\right]$.
If $\deg_{A}\left(Q\right)>0$, then the valuation/pole-order test
shows that $T\left[G\right]$ has a finite pole in the $A$-line,
contradicting $T\left[G\right]=\left(n+1\right)A\in\widetilde{\mathcal{K}}\left[A\right]$.
Consequently, $\deg_{A}\left(Q\right)=0$ and $G\in\widetilde{\mathcal{K}}\left[A\right]$. 
\end{singlespace}

Now set
\[
G={\displaystyle \sum_{j=0}^{d}\alpha_{j}A^{j}}
\]
with all $\alpha_{j}\in\widetilde{\mathcal{K}}$. If $d\geq2$, then
comparing the $A^{d}$-coefficients in $T\left[G\right]=\left(n+1\right)A$
yields
\[
T\left[\alpha_{d}\right]=0,
\]
implying $\alpha_{d}\in\mathbb{C}$. Furthermore, comparing the $A^{d-1}$-coefficients
reveals
\[
T\left[\alpha_{d-1}\right]+d\alpha_{d}\left(nH_{n-1}\right)=0,
\]
which is impossible by Assertion $\left(\textrm{II}_{n}\right)$ unless
$\alpha_{d}=0$, forcing $d\leq1$. With $G=\alpha_{1}A+\alpha_{0}$,
a comparison of the $A^{1}$-coefficients gives 
\[
T\left[\alpha_{1}\right]=n+1
\]
\[
\Longleftrightarrow\alpha_{1}^{\prime}\left(t\right)=\tfrac{n+1}{t},
\]
which has no solution $\alpha_{1}\in\widetilde{\mathcal{K}}$ because
every element of $\widetilde{\mathcal{K}}$ is meromorphic at $t=0$
with a Laurent expansion free of logarithms. This contradiction validates
Assertion $\left(\textrm{II}_{n+1}\right)$.\medskip{}

\begin{singlespace}
Through the above Steps, we have shown by induction that both Assertions
$\left(\textrm{I}_{n}\right)$ and $\left(\textrm{II}_{n}\right)$
hold for all $n\geq1$. In particular,
\[
\mathrm{Const}\left(\mathcal{K}_{n},T\right)=\mathbb{C}
\]
for all $n\geq0$. $\blacklozenge$\medskip{}
\medskip{}

\end{singlespace}

\begin{singlespace}
\noindent\textbf{Fact 3:} $H_{1}\notin\mathcal{K}\left(E\right)$.

\noindent Recalling that $E=e^{-t}$ is transcendental over $\mathcal{K}$
and $T[H_{1}]=H_{0}=1-E\in\mathcal{K}[E]$, let us assume (for purposes
of contradiction) that $H_{1}\in\mathcal{K}(E)$. This implies there
exists $G\in\mathcal{K}\left(E\right)$ such that
\begin{equation}
T\left[G\right]=1-E.
\end{equation}

\end{singlespace}

\begin{singlespace}
Note that we work in an algebraic differential extension where $T$
extends and all denominators split. (In algebraic extensions of $\left(\mathcal{K},T\right)$,
the field of constants enlarges only algebraically, which is sufficient
because $Ce^{-t}$ is algebraic over $\mathbb{C}(t)$ only if $C=0$.)
First, write the partial-fraction expansion
\[
G=U\left(E\right)+\sum_{\rho}\sum_{j=1}^{\mu_{\rho}}\frac{a_{\rho,j}}{\left(E-\rho\right)^{j}},
\]
where $U\left(E\right)\in\mathcal{K}\left[E\right]$ and $a_{\rho,j},\rho\in\overline{\mathcal{K}}$.
Next, fix a pole $\rho$ and let $\mu_{\rho}=\mathrm{ord}_{E=\rho}\left(G\right)\ge1$,
so that $a_{\rho,\mu_{\rho}}\neq0$ and $a_{\rho,\mu_{\rho}+1}=0$.
Given that $T\left[E\right]=-tE$, we find
\[
T\left[\left(E-\rho\right)^{-j}\right]=-j\left(E-\rho\right)^{-j-1}T\left[E-\rho\right]
\]
\[
=j\left(E-\rho\right)^{-j-1}\left[t\left(E-\rho\right)+\left(t\rho+T\left[\rho\right]\right)\right].
\]
It then follows that the coefficient, $\left(E-\rho\right)^{-\left(\mu_{\rho}+1\right)}$,
of the highest pole in $T\left[G\right]$ equals\linebreak{}
$\mu_{\rho}a_{\rho,\mu_{\rho}}\left(t\rho+T\left[\rho\right]\right)$
because the potential contribution from $T\left[a_{\rho,\mu_{\rho}+1}\right]$
vanishes. Since $T\left[G\right]=1-E\in\mathcal{K}\left[E\right]$
has no poles in $E$, we obtain $t\rho+T\left[\rho\right]=0$ for
each pole $\rho$ of $G$. Solving $T\left[\rho\right]=-t\rho$ (i.e.,
$\rho'=-\rho$) then gives $\rho=Ce^{-t}$ (with $C$ a constant in
the extension), which forces $\rho=0$. Thus, all poles of $G$ (if
any) occur at $E=0$, and we can write a Laurent expansion
\[
G=\sum_{k=-\mu}^{N}b_{k}\left(t\right)E^{k},
\]
where $\mu\geq0$ and $b_{-\mu}\left(t\right)\neq0$ if $\mu>0$.
Since $E$ is transcendental over $\mathcal{K}$, we see that the
monomials $\left\{ E^{k}:k\in\mathbb{Z}\right\} $ are $\mathcal{K}$-linearly
independent, making the Laurent expansion unique.

Applying $T$ to this expansion gives
\[
T\left[G\right]=\sum_{k=-\mu}^{N}\left(T\left[b_{k}\right]-ktb_{k}\right)E^{k}.
\]
Then, since $T\left[G\right]\in\mathcal{K}\left[E\right]$, the coefficients
at $E^{k}$ must vanish for $k<0$; that is:
\[
T\left[b_{k}\right]-ktb_{k}=0.
\]
Equivalently, $b_{k}^{\prime}=kb_{k}$ (because $t\in\mathcal{K}^{\times}$),
and the only rational solution of this first-order differential equation
with $k\in\mathbb{C}^{\times}$ is $b_{k}=0$. (Indeed, if $b_{k}\in\mathbb{C}\left(t\right)^{\times}$,
then $b^{\prime}/b$ is rational with only simple poles at the zeros/poles
of $b$; and if $b^{\prime}/b=k$ is constant, then there are no poles,
implying $b\in\mathbb{C}^{\times}$ and $b^{\prime}=0\neq kb$.) The
fact that $b_{k}=0$ for all $k<0$ forces $\mu=0$, and thus $G\in\mathcal{K}\left[E\right]$.

Now write
\[
G=\sum_{k=0}^{N}b_{k}\left(t\right)E^{k},
\]
so that
\begin{equation}
T\left[G\right]=\sum_{k=0}^{N}\left(T\left[b_{k}\right]-ktb_{k}\right)E^{k}.
\end{equation}

Since the $\left\{ E^{k}\right\} _{k\geq0}$ are $\mathcal{K}$-linearly
independent, comparing (18) with (17) yields
\[
T\left[b_{0}\right]=1.
\]
However, $T\left[b_{0}\right]=1$ implies $b_{0}^{\prime}=1/t$, which
has no solution in $\mathbb{C}\left(t\right)$ (by the residue test
at $t=0$). This contradiction shows that the posited $G\in\mathcal{K}\left(E\right)$
cannot exist, and so $H_{1}\notin\mathcal{K}\left(E\right)$. $\blacklozenge$\medskip{}
\medskip{}

\end{singlespace}

\begin{singlespace}
\noindent\textbf{Fact 4:} $E$ is transcendental over $\mathcal{K}\left(H_{1},H_{2},\ldots,H_{m}\right)$
for any $m\in\mathbb{Z}_{\geq1}$.

\noindent We work with the fixed derivation $T$ on the tower
\[
\mathcal{K}_{0}\subset\mathcal{K}_{1}\subset\mathcal{K}_{2}\subset\cdots\subset\mathcal{K}_{m+1},
\]
recalling that: (i) by Fact 2, $\textrm{Const}\left(\mathcal{K}_{n},T\right)=\mathbb{C}$
for all $n\geq0$; (ii) $E$ is an exponential of an integral (i.e.,
$T\left[E\right]/E=-t\in\mathcal{K}_{0}$); and (iii) for $n\geq1$,
each $H_{n}$ is adjoined as an antiderivative (i.e., $T\left[H_{n}\right]=nH_{n-1}\in\mathcal{K}_{n-1}$).
\end{singlespace}

\begin{singlespace}
Assume (for purposes of contradiction) that $E$ is algebraic over
$\mathcal{K}_{m+1}$ for some $m\in\mathbb{Z}_{\geq1}$, and let
\[
P\left(Z\right)={\displaystyle \sum_{j=0}^{d}}\alpha_{j}\left(H_{1},H_{2},\ldots,H_{m}\right)Z^{j}\in\mathcal{K}\left(H_{1},H_{2},\ldots,H_{m}\right)\left[Z\right]
\]
denote the monic minimal polynomial of $E$ over $\mathcal{K}\left(H_{1},H_{2},\ldots,H_{m}\right)$.
Further, choose a common denominator $J\left(H_{1},H_{2},\ldots,H_{m}\right)\in\mathcal{K}\left[H_{1},H_{2},\ldots,H_{m}\right]$
for the coefficients $\alpha_{j}$, set $\widetilde{\alpha}_{j}=J\alpha_{j}\in\mathcal{K}\left[H_{1},H_{2},\ldots,H_{m}\right]$,
and define
\[
Q\left(\mathbf{W},Z\right)=\sum_{j=0}^{d}\widetilde{\alpha}_{j}(\mathbf{W})Z^{j}\in\mathcal{K}\left[\mathbf{W},Z\right]
\]
for $\mathbf{W}=\left(W_{1},W_{2},\ldots,W_{m}\right)$. Since $P$
is monic in $Z$, we have $Q\neq0$ as a polynomial in the indeterminates
$(\mathbf{W},Z)$.

Given that $E$ is transcendental over $\mathcal{K}$ by Fact 1, the
specialization map
\[
\mathrm{ev}_{E}:\mathcal{K}\left[\mathbf{W},Z\right]\longrightarrow\mathcal{K}\left(E\right)\mathbf{\left[W\right]}
\]
for $R\left(\mathbf{W},Z\right)\mapsto R\left(\mathbf{W},E\right)$
is injective. (Indeed, if $R\left(\mathbf{W},Z\right)=0$ with $R\left(\mathbf{W},Z\right)={\textstyle \sum}_{\iota}h_{\iota}\left(\mathbf{W}\right)Z^{\iota}$
and $h_{\iota}\in\mathcal{K}\left[\mathbf{W}\right]$, then for each
monomial $M\left(\mathbf{W}\right)$ the coefficient $\sum_{\nu}\beta_{\nu,M}E^{\nu}$
must vanish in $\mathcal{K}(E)$. This is because the transcendence
of $E$ forces all $\beta_{\nu,M}=0$, and thus all $h_{\iota}=0$,
yielding $R=0$.) Consequently,
\[
R\left(\mathbf{W}\right)=Q\left(\mathbf{W},E\right)\in\mathcal{K}\left(E\right)\mathbf{\left[W\right]}
\]
is a nonzero polynomial, and evaluating at $\mathbf{W}=\mathbf{H}=\left(H_{1},H_{2},\ldots,H_{m}\right)$
gives
\[
R\left(\mathbf{H}\right)=Q\left(\mathbf{H},E\right)
\]
\[
=J\left(\mathbf{H}\right)P\left(E\right)
\]
\[
=0;
\]
that is, a genuine algebraic relation among the $H_{1},H_{2},\ldots,H_{m}$
over the base field $\mathcal{K}(E)$:
\[
R\left(H_{1},H_{2},\ldots,H_{m}\right)=0
\]
with $R\in\mathcal{K}\left(E\right)\left[W_{1},W_{2},\ldots,W_{m}\right],\:R\neq0$.

Now apply the Kolchin-Ostrowski theorem for towers of antiderivatives
(e.g., Theorem 4.6 of Srinivasan, 2008) over the base differential
field $\left(\mathcal{K}\left(E\right),T\right)$ whose constants
are $\mathbb{C}$. Since $H_{1},H_{2},\ldots,H_{m}$ are obtained
by iterated adjunction of $T$-antiderivatives and are algebraically
dependent over $\mathcal{K}\left(E\right)$, there exist constants
$c_{1},c_{2},\ldots,c_{m}\in\mathbb{C}$, at least one of which is
nonzero, such that
\[
S_{1}=c_{1}H_{1}+c_{2}H_{2}+\cdots+c_{m}H_{m}\in\mathcal{K}\left(E\right).
\]
Letting $k^{*}=\max\left\{ k:c_{k}\neq0\right\} $ and applying $T^{k^{*}-1}$
then yields
\[
T^{k^{*}-1}\left[S_{1}\right]=k^{*}!c_{k^{*}}H_{1}+V_{1}\in\mathcal{K}\left(E\right)
\]
for some $V_{1}\in\mathcal{K}\left(E\right)$ (since $H_{0}=1-E\in\mathcal{K}\left(E\right)$),
forcing $H_{1}\in\mathcal{K}\left(E\right)$, which contradicts Fact
3. Therefore, $E$ cannot be algebraic over $\mathcal{K}\left(H_{1},H_{2},\ldots,H_{m}\right)$,
and so must be transcendental over that field. $\blacklozenge$\medskip{}
\medskip{}

Given the above facts, we assume (for purposes of contradiction) that
the $H_{1},H_{2},\ldots,H_{m}$ satisfy a nontrivial algebraic relation
over $\mathcal{K}\left(E\right)$. Recalling Fact 2, and applying
the Kolchin-Ostrowski theorem for towers of antiderivatives (e.g.,
Theorem 4.6 of Srinivasan, 2008 -- applicable for any derivation,
in particular $T$) over $\left(\mathcal{K}\left(E\right),T\right)$
implies the existence of constants $\widetilde{c}_{1},\widetilde{c}_{2},\ldots,\widetilde{c}_{m}\in\mathbb{C}$,
at least one of which is nonzero, such that
\[
S_{2}=\widetilde{c}_{1}H_{1}+\widetilde{c}_{2}H_{2}+\cdots+\widetilde{c}_{m}H_{m}\in\mathcal{K}\left(E\right).
\]
Letting $\ell^{*}=\max\left\{ \ell:\widetilde{c}_{\ell}\neq0\right\} $
and applying $T^{\ell^{*}-1}$ then gives
\[
T^{\ell^{*}-1}\left[S_{2}\right]=\ell^{*}!\widetilde{c}_{\ell^{*}}H_{1}+V_{2}\in\mathcal{K}\left(E\right)
\]
for some $V_{2}\in\mathcal{K}\left(E\right)$ (since $H_{0}=1-E\in\mathcal{K}\left(E\right)$),
which means $H_{1}\in\mathcal{K}\left(E\right)$, contradicting Fact
3. Therefore, the $H_{1},H_{2},\ldots,H_{m}$ must be algebraically
independent over $\mathcal{K}\left(E\right)$.

Now assume (again for contradiction) that there exists a nontrivial
polynomial relation over $\mathcal{K}$ among the $m+1$ functions
$E,H_{1},H_{2},\ldots,H_{m}$. Since viewing this as a polynomial
in $E$ with coefficients in $\mathcal{K}\left[H_{1},H_{2},\ldots,H_{m}\right]$
would make $E$ algebraic over $\mathcal{K}\left(H_{1},H_{2},\ldots,H_{m}\right)$
-- thereby contradicting Fact 4 -- we conclude that the $E,H_{1},H_{2},\ldots,H_{m}$
must be algebraically independent over $\mathcal{K}$.

Finally, recall from (14) and (15) that $E=1-tF_{0}$ and $H_{n}=tF_{n},\:n\geq1$,
respectively. This implies that the two $\left(m+1\right)$-tuples,
$\left(E,H_{1},H_{2},\ldots,H_{m}\right)$ and $\left(F_{0},F_{1},\ldots,F_{m}\right)$,
generate the same subfield and thus possess the same transcendence
degree. In particular, if the $F_{0},F_{1},\ldots,F_{m}$ satisfied
a nontrivial algebraic relation over $\mathcal{K}$, then substituting
$F_{0}=\left(1-E\right)/t$ and $F_{n}=H_{n}/t$ would create a nontrivial
algebraic relation among the $E,H_{1},H_{2},\ldots,H_{m}$, which
we just have ruled out. Therefore, the $F_{0},F_{1},\ldots,F_{m}$
must be algebraically independent over $\mathbb{C}\left(t\right)$,
the desired conclusion. $\blacksquare$ \medskip{}

Theorems 1 and 2 provide our main results, with special cases involving
$\gamma$ and $\delta$ addressed in ancillary corollaries.\medskip{}
\medskip{}

\end{singlespace}

\begin{singlespace}
\noindent\textbf{Theorem 1:} The set of values $\left\{ \eta^{\left(n\right)}\right\} _{n\geq0}$
are algebraically independent over $\overline{\mathbb{Q}}$, implying
\[
\eta^{\left(n\right)}\in\mathbb{R}\setminus\overline{\mathbb{Q}}
\]
for all $n\in\mathbb{Z}_{\geq0}$.\medskip{}
\medskip{}
\medskip{}

\noindent\textbf{Proof:} For any $m\in\mathbb{Z}_{\geq1}$, define
\[
Y_{n}\left(t\right)=F_{n-1}\left(t\right)
\]
for $n\in\left\{ 1,2,\ldots,m+1\right\} $ and
\[
Y_{m+2}\left(t\right)=1.
\]
Equations (12) and (13) then generate the homogeneous linear system
\begin{equation}
\boldsymbol{Y}^{\prime}\left(t\right)=\mathbf{A}\left(t\right)\boldsymbol{Y}\left(t\right),
\end{equation}
where $\boldsymbol{Y}=\left[Y_{1},Y_{2},\ldots,Y_{m+2}\right]^{\textrm{T}}$
and $\mathbf{A}\left(t\right)\in M_{m+2}\left(\mathbb{Q}\left(t\right)\right)$
is the matrix
\[
\mathbf{A}\left(t\right)=\left[\begin{array}{ccccccccc}
-1-\dfrac{1}{t} & 0 & 0 & 0 & \cdots & 0 & 0 & 0 & \dfrac{1}{t}\\
\dfrac{1}{t} & -\dfrac{1}{t} & 0 & 0 & \cdots & 0 & 0 & 0 & 0\\
0 & \dfrac{2}{t} & -\dfrac{1}{t} & 0 & \cdots & 0 & 0 & 0 & 0\\
0 & 0 & \dfrac{3}{t} & -\dfrac{1}{t} & \cdots & 0 & 0 & 0 & 0\\
\vdots & \vdots & \vdots & \vdots & \ddots & \vdots & \vdots & \vdots & \vdots\\
0 & 0 & 0 & 0 & \cdots & -\dfrac{1}{t} & 0 & 0 & 0\\
0 & 0 & 0 & 0 & \cdots & \dfrac{m-1}{t} & -\dfrac{1}{t} & 0 & 0\\
0 & 0 & 0 & 0 & \cdots & 0 & \dfrac{m}{t} & -\dfrac{1}{t} & 0\\
0 & 0 & 0 & 0 & \cdots & 0 & 0 & 0 & 0
\end{array}\right]
\]
for $t\in\mathbb{R}^{+}$. We note that the only singularity of the
system occurs at $t=0$.
\end{singlespace}

\begin{singlespace}
From Lemma 1, we know that the $E$-functions $F_{0}\left(t\right),F_{1}\left(t\right),\ldots,F_{m}\left(t\right)$
are algebraically independent over $\mathbb{C}\left(t\right)$. Since
$\overline{\mathbb{Q}}\left(t\right)\subset\mathbb{C}\left(t\right)$,
the functions must be algebraically independent over $\overline{\mathbb{Q}}\left(t\right)$
as well. Therefore, the vector $\left[F_{0}\left(t\right),F_{1}\left(t\right),\ldots,F_{m}\left(t\right),1\right]^{\textrm{T}}$
consists of $m+2$ $E$-functions (including $\widetilde{F}\left(t\right)=1$)
that solve a linear system with coefficients in $\mathbb{Q}\left(t\right)$,
and whose first $m+1$ components are algebraically independent.

By applying the (refined) Siegel-Shidlovski theorem (see Beukers,
2006) at the algebraic point $t=1$ (which is not a singularity of
the system), we know that every algebraic relation over $\overline{\mathbb{Q}}$
among the values
\[
\left[F_{0}\left(1\right),F_{1}\left(1\right),\ldots,F_{m}\left(1\right)\right]^{\textrm{T}}=\left[\eta^{\left(0\right)},\eta^{\left(1\right)},\ldots,\eta^{\left(m\right)}\right]^{\textrm{T}}
\]
must arise from an algebraic relation over $\overline{\mathbb{Q}}\left(t\right)$
among the functions $F_{0}\left(t\right),F_{1}\left(t\right),\ldots,F_{m}\left(t\right)$.\footnote{The purpose of including $Y_{m+2}\left(t\right)=1$ in (19) is simply
to embed the $m+1$ functions $F_{0}\left(t\right),F_{1}\left(t\right),\ldots,F_{m}\left(t\right)$
into a single linear system. The Siegel-Shidlovski theorem then allows
us to deduce the algebraic independence of the corresponding $m+1$
values at $t=1$.} Since these functions are algebraically independent over $\overline{\mathbb{Q}}\left(t\right)$,
we can conclude that the set of values\linebreak{}
$\eta^{\left(0\right)},\eta^{\left(1\right)},\ldots,\eta^{\left(m\right)}$
are algebraically independent over $\overline{\mathbb{Q}}$, and thus
that the individual $\eta^{\left(n\right)}$ are transcendental. Taking
$m$ arbitrarily large provides the final result. $\blacksquare$\medskip{}

\end{singlespace}

\begin{singlespace}
\noindent\textbf{Corollary 1:}
\[
\eta=\gamma+\dfrac{\delta}{e}\in\mathbb{R}\setminus\overline{\mathbb{Q}},
\]
where $\gamma$ and $\delta$ denote the Euler-Mascheroni and Euler-Gompertz
constants, respectively.\medskip{}
\medskip{}

\noindent\textbf{Proof:} This result follows as a special case of
Theorem 1 for $n=1$. $\blacksquare$
\end{singlespace}

\begin{singlespace}
\medskip{}

\end{singlespace}

\begin{singlespace}
\noindent\textbf{Theorem 2:} For all $n\in\mathbb{Z}_{\geq1}$:
\[
\left(\textrm{i}\right)\:\left(\gamma^{\left(n\right)}\in\mathbb{R}\setminus\overline{\mathbb{Q}}\right)\vee\left(\dfrac{\delta^{\left(n\right)}}{e}\in\mathbb{R}\setminus\overline{\mathbb{Q}}\right);\textrm{and}
\]
\[
\left(\textrm{ii}\right)\:\left(\gamma^{\left(n\right)}\in\mathbb{R}\setminus\overline{\mathbb{Q}}\right)\vee\left(\delta^{\left(n\right)}\in\mathbb{R}\setminus\overline{\mathbb{Q}}\right).
\]
\medskip{}

\noindent\textbf{Proof:\medskip{}
}

\noindent (i) Theorem 1 shows that $\eta^{\left(n\right)}$ is transcendental
for all $n\in\mathbb{Z}_{\geq0}$. Moreover, it is clear from equation
(6) that if $\eta^{\left(n\right)}$ is transcendental for any given
$n\geq0$, then the two numbers $\gamma^{\left(n\right)}$ and $\delta^{\left(n\right)}/e$
cannot both be algebraic. (We omit the case of $n=0$ from the statement
of the theorem because it is trivially obvious that $\delta^{\left(0\right)}/e=-e^{-1}$
is transcendental.)\medskip{}

\noindent (ii) From Theorem 1, we know that $\eta^{\left(n\right)}=\gamma^{\left(n\right)}+\delta^{\left(n\right)}/e$
is algebraically independent of $\eta^{\left(0\right)}=1-e^{-1}$
for all $n\in\mathbb{Z}_{\geq1}$. Therefore, if we assume (for purposes
of contradiction) that both $\gamma^{\left(n\right)}$ and $\delta^{\left(n\right)}$
are algebraic, then we can write 
\[
\eta^{\left(n\right)}=\gamma^{\left(n\right)}+\delta^{\left(n\right)}\left(1-\eta^{\left(0\right)}\right)
\]
\[
\Longleftrightarrow\eta^{\left(n\right)}-\left(\gamma^{\left(n\right)}+\delta^{\left(n\right)}\right)+\delta^{\left(n\right)}\eta^{\left(0\right)}=0
\]
\[
\Longleftrightarrow a\eta^{\left(n\right)}+b\eta^{\left(0\right)}+c=0,
\]
where $a,b,c\in\overline{\mathbb{Q}}$. This contradicts the algebraic
independence of $\eta^{\left(n\right)}$ and $\eta^{\left(0\right)}$,
forcing either $\gamma^{\left(n\right)}$ or $\delta^{\left(n\right)}$
to be transcendental. $\blacksquare$\medskip{}
\medskip{}

\noindent\textbf{Corollary 2:}
\[
\left(\textrm{i}\right)\:\left(\gamma\in\mathbb{R}\setminus\overline{\mathbb{Q}}\right)\vee\left(\dfrac{\delta}{e}\in\mathbb{R}\setminus\overline{\mathbb{Q}}\right)\textrm{and}
\]
\[
\left(\textrm{ii}\right)\:\left(\gamma\in\mathbb{R}\setminus\overline{\mathbb{Q}}\right)\vee\left(\delta\in\mathbb{R}\setminus\overline{\mathbb{Q}}\right),
\]
where $\gamma$ and $\delta$ denote the Euler-Mascheroni and Euler-Gompertz
constants, respectively.\medskip{}
\medskip{}

\noindent\textbf{Proof:} These results follow as special cases of
the corresponding parts of Theorem 2 for $n=1$. As noted in the Introduction,
(ii) is already known from Rivoal (2012). $\blacksquare$\medskip{}
\medskip{}

\end{singlespace}
\begin{singlespace}

\section{Conclusion}
\end{singlespace}

\begin{singlespace}
\noindent In the present article, we defined sequences of generalized
Euler-Mascheroni, Euler-Gompertz, and Eta constants, denoted by $\gamma^{\left(n\right)}$,
$\delta^{\left(n\right)}$, and $\eta^{\left(n\right)}$, respectively.
After summarizing the basic properties of these sequences, we provided
two new results. First, Theorem 1 proved the algebraic independence
of the set $\left\{ \eta^{\left(n\right)}\right\} _{n\geq0}$, and
thus the transcendence of $\eta^{\left(n\right)}$ for all $n\in\mathbb{Z}_{\geq0}$.
This includes the relatively prominent value $\eta=\gamma+\delta/e$
introduced by Hardy (1949). Second, Theorem 2 showed that the transcendence
result of Theorem 1 implies the disjunctive transcendence of both
pairs $\left(\gamma^{\left(n\right)},\delta^{\left(n\right)}/e\right)$
and $\left(\gamma^{\left(n\right)},\delta^{\left(n\right)}\right)$
for all $n\in\mathbb{Z}_{\geq1}$. For the special case of $\left(\gamma,\delta/e\right)$,
this complements the disjunctive transcendence of $\left(\gamma,\delta\right)$
(proved by Rivoal, 2012) and also draws attention to the fact that,
under the probabilistic interpretation of Hardy's equation, the constant
$\delta/e=\textrm{E}_{X}\left[X^{\left(-\right)}\right]=-\textrm{E}_{X}\left[\min\left\{ X,0\right\} \right]$
(for $X\sim\textrm{Gumbel}\left(0,1\right)$) may be a more natural
companion of $\gamma=\textrm{E}_{X}\left[X\right]$ than is $\delta=-\textrm{E}_{X}\left[X\mid X\leq0\right]$.
\end{singlespace}

\begin{singlespace}
Although the primary focus of the current work has been the transcendence
of the three constant sequences discussed above, it is important to
mention that Theorem 1 could be restated to prove the algebraic independence
of the set of specialized values $\left\{ F_{n}\left(t\right)\right\} _{n\geq0}$
(and thus the transcendence of the individual values $F_{n}\left(t\right)$)
for \emph{any} $t\in\overline{\mathbb{Q}}\setminus0$. One particularly
interesting case is that of $t=-1$, for which the sequence of constants
\[
F_{n}\left(-1\right)=n!{\displaystyle \sum_{k=1}^{\infty}\dfrac{1}{k^{n}\cdot k!}}
\]
includes
\begin{equation}
\sum_{k=1}^{\infty}\dfrac{1}{k\cdot k!}=-\left(\gamma+\dfrac{\delta^{*}}{e}\right),
\end{equation}
where
\[
\delta^{*}\equiv-e\textrm{Ei}\left(1\right)
\]
\[
=-5.151464\ldots.
\]
Powers (2025) used (20) to motivate $\delta^{*}$ as an ``alternating
analogue'' of $\delta$, with a relationship comparable to that between
$\ln\left(\tfrac{4}{\pi}\right)=0.241564\ldots$ and $\gamma$ (see
Sondow, 2005).\medskip{}
\medskip{}

\end{singlespace}

\end{document}